\documentclass[11pt,fleqn,a4paper]{article}
\usepackage{amsmath,amssymb,amsthm,hyperref} 

\newcommand{\p}{{\partial}}

\newcommand{\diag}{\mathop{\rm diag}\nolimits}

\newtheorem{theorem}{Theorem}

\newtheorem{corollary}{Corollary}

{\theoremstyle{definition}

\newtheorem{remark}{Remark}
\newtheorem{example}{Example}
\newtheorem*{remark*}{Remark}

}

\newcommand{\todo}[1][\null]{\ensuremath{\clubsuit}}

\newcommand{\noprint}[1]{}

\begin{document}

\begin{flushleft}
\LARGE \bf
Lie symmetries of systems of second-order linear ordinary differential equations with constant coefficients
\end{flushleft}

\begin{flushleft}
Vyacheslav M.~BOYKO~$^\dag$, Roman O.~POPOVYCH~$^{\dag\ddag}$ and Nataliya M. SHAPOVAL~$^{\S}$
\end{flushleft}

\noindent $^\dag$~Institute of Mathematics of NAS of Ukraine, 3
Tereshchenkivs'ka Str., Kyiv-4, 01601 Ukraine\\
$\phantom{^\dag}$~E-mail: boyko@imath.kiev.ua, rop@imath.kiev.ua

\medskip

\noindent
$^\ddag$~Wolfgang Pauli Institut, Universit\"at Wien, Nordbergstra{\ss}e 15, A-1090 Wien, Austria

\medskip

\noindent
$^\S$~Faculty of Mechanics and Mathematics, National Taras Shevchenko University of Kyiv,\\
\hphantom{${}^\S$} 2 Academician Glushkov Ave., 03127 Kyiv, Ukraine\\
$\phantom{^\S}$~E-mail: natalya.shapoval@gmail.com

{\vspace{7mm}\par\noindent\hspace*{5mm}\parbox{150mm}{\small
Lie symmetries of systems of second-order linear ordinary differential equations with commuting constant matrix coefficients are exhaustively described
over both the complex and real fields.
The exact lower and upper bounds for the dimensions of the maximal Lie invariance algebras possessed by such systems
are obtained using an effective algebraic approach.
}\par\vspace{7mm}}

\noindent
{\it Key words:} system of second-order ordinary differential equation; Lie symmetry; matrix equation

\bigskip

\noindent
{\it \it 2010 Mathematics Subject Classification:} 34A30; 34C14

\section{Introduction}

The problem on possible dimensions of the maximal Lie invariance algebras of differential equations from a fixed class has a long history.
Already S.~Lie obtained exhaustive results concerning maximal dimensions of such algebras for ordinary differential equations (ODEs) of any fixed order \cite[S.~294--301]{Lie1893}.
Namely, he first proved that any first-order ODE possesses an infinite-dimensional Lie invariance algebra,
the dimension of the maximal Lie invariance algebra of any second-order ODE (resp.\ any $m$th order ODE for $m\geqslant 3$) is at most eight (resp.\ not greater than $m+4$),
and these bounds are exact.
Later these results were repeatedly reinvented, see e.g.~\cite{Gonzalez-Gascon&Gonzalez-Lopez1983}.

Analogous results for systems of ODEs are much less known. We discuss some of them, which are relevant to the subject of the present paper.
Thus, according to the remarkable lecture notes by Markus \cite[pp.~68--69, Theorem 44]{Markus1960},
any system of second-order ODEs
\begin{equation}\label{second_order_system}
\ddot{\boldsymbol{x}} = \boldsymbol{f}(t,\boldsymbol{x},\dot{\boldsymbol{x}}),
\end{equation}
where
$\boldsymbol{x}(t)=(x^1(t),\dots,x^n(t))^{\rm T}$,
$\dot{\boldsymbol{x}} =d\boldsymbol{x}/dt$,
$\ddot{\boldsymbol{x}}=d\dot{\boldsymbol{x}}/dt$,
possesses the maximal Lie invariance algebra of dimensions not greater than $(n+2)^2-1$.
This result was later reproved in \cite[Sections~4 and~5]{Gonzalez-Gascon&Gonzalez-Lopez1983}.
It was also shown therein that the maximal dimension $(n+2)^2-1=n^2+4n+3$ is reached
for systems reduced by point transformations to the simplest system
\begin{equation}\label{x''=0}
\ddot{\boldsymbol{x}}=\boldsymbol{0}.
\end{equation}
The maximal Lie invariance algebra~$\mathfrak g^0$ of the system \eqref{x''=0} is generated by the vector fields
\begin{gather*}
\p_t,\quad \p_{x^a}, \quad t\p_t, \quad x^a\p_t,\quad t\p_{x^a},
\quad x^a\p_{x^b}, \quad tx^a\p_t+x^ax^c\p_{x^c}, \quad t^2\p_t+tx^c\p_{x^c},
\end{gather*}
and is isomorphic to the Lie algebra
$\mathrm{sl}(n+2,{\mathbb C})$ (resp.\ $\mathrm{sl}(n+2,{\mathbb R})$)
for the complex (resp.\ real) case, see e.g.~\cite{Gonzalez-Lopez1988}.
Here and in what follows the indices $a$, $b$, $c$ run from 1 to $n$, i.e.\ $a,b,c=1,\dots,n$,
and we use the summation convention for repeated indices.
It can be checked that the system~\eqref{x''=0} is invariant with respect to
the general projective group of ${\mathbb C}^{n+1}$ (resp.\ ${\mathbb R}^{n+1}$)
consisting of the point transformations \mbox{\cite[S.~554]{Lie1888}}
\[
\tilde x^i=\frac{\alpha_{i0}x^0+\dots+\alpha_{in}x^n+\alpha_{i,n+1}}{\alpha_{n+1,0}x^0+\dots+\alpha_{n+1,n}x^n+\alpha_{n+1,n+1}},\qquad i=0,\dots,n,
\]
where $\alpha_{00}$, $\alpha_{01}$, \dots, $\alpha_{n+1,n+1}$ are homogeneous group parameters and $x^0=t$,
and $\mathfrak g^0$ is the Lie algebra associated with this group.
In fact, the number of essential group parameters is $(n+2)^2-1$ as, supposing a homogeneous group parameter nonzero,
we can set it to be equal to~$1$, simultaneously dividing of the numerator and the denominator in the expression for each~$\tilde x^i$ by this parameter
and then assigning parameter ratios as new parameters.

Fels proved \cite{Fels1993,Fels1995} that up to point equivalence the system \eqref{x''=0}
is a unique system of the form~\eqref{second_order_system} which admits an $(n^2+4n+3)$-dimensional Lie invariance algebra,
which was earlier known only for linear systems~\cite{Gonzalez-Lopez1988}.
Recently linearization criteria for systems from the class~\eqref{second_order_system} have been independently investigated in \cite{Aminova2010,Bagderina2010,Merker2006}.
See also the references therein.
The maximal dimension of the maximal Lie invariance algebras for normal systems of $m$th order ODEs
was estimated in \cite{Gonzalez-Gascon&Gonzalez-Lopez1983,Gonzalez-Gascon&Gonzalez-Lopez1985,Gonzalez-Gascon&Gonzalez-Lopez1988}
for an arbitrary $m\geqslant3$,
and it is known for the case $m=3$ that up to point equivalence the system $\dddot{\boldsymbol{x}}=0$ is a unique system
for which the dimension of the maximal Lie invariance algebra reaches the maximal value $n^2+3n+3$ for such systems~\cite{Fels1993}.

In a series of recent papers the study of Lie symmetries of systems of $n$ ($n\geqslant 2$) linear second-order ODEs with constant coefficients was recovered.
Namely, the cases $n=2$ and $n=3$ were considered in~\cite{WafoSoh2010}.
In~\cite{Campoamor2011} Campoamor-Stursberg corrected results of~\cite{WafoSoh2010}
(see also comments concerning~\cite{WafoSoh2010} in~\cite{Meleshko2011}),
and studied the case $n=4$ as well as systems associated with diagonal matrices without restrictions on~$n$.
Certain results on the dimensions of the maximal Lie invariance algebras of such systems
in the case of arbitrary $n$ and matrices of the general Jordan form were obtained in~\cite{Campoamor2012}.

The study of symmetry properties of systems from the class~\eqref{second_order_system} is required by numerous applications in mechanics, gravity, etc.
Unfortunately, there are no general results on Lie symmetries of these systems.
This is why even linear systems with constant coefficients are good objects for a preliminary investigation
in spite of the well-known simple algorithm for the construction of their general solutions.
Group classification of linear systems with constant coefficients gives examples what
Lie algebras of vector fields are admitted by systems from the class~\eqref{second_order_system} as their maximal Lie invariance algebras
and what dimensions of these algebras are possible.
Note that the above knowledge is important for the problem on linearization of systems from the class~\eqref{second_order_system}.

The purpose of the present paper, which is inspired by the papers~\cite{Campoamor2011,Campoamor2012,Meleshko2011,WafoSoh2010},
is to exhaustively describe Lie symmetries of systems of second-order linear ordinary differential equations with commuting constant matrix coefficients
over both the complex and real fields.
We essentially enhance and generalize the results of~\cite{Campoamor2011,Campoamor2012,Meleshko2011,WafoSoh2010}
using a simple but effective algebraic approach.
In particular, we explicitly describe the maximal Lie invariance algebras possessed by systems under consideration
with no restrictions on the number of equations and the form of system coefficients and
derive the exact lower and upper bounds for the dimensions of these algebras.

\section{Main result}

Following \cite{Campoamor2011, Campoamor2012,Meleshko2011,WafoSoh2010}, we consider systems of linear second-order ODEs of the normal form
\begin{gather}\label{GeneralLinearSystemODEs}
\ddot{\boldsymbol{x}}=A\dot{\boldsymbol{x}}+B\boldsymbol{x}+{\boldsymbol{C}}(t)
\end{gather}
over the complex field.
Here ${\boldsymbol{C}}(t)$ is a smooth $n$-component vector-function of~$t$,
and $A$ and $B$ are \emph{commuting} constant complex-valued matrices of dimension $n\times n$, $n\geqslant 2$.
Note that the choice of the underlying field ($\mathbb C$ or $\mathbb R$) is not principal.
We choose the complex field in order to make the presentation clearer.

It is commonly known (see, e.g., \cite{Gonzalez-Lopez1988})
that the change of dependent variables $\boldsymbol{x}=\exp(\frac{1}{2}At)\boldsymbol{y}+{\boldsymbol{x}}_{\text{p}}(t)$, where ${\boldsymbol{x}}_{\text{p}}(t)$ is
a~particular solution of~\eqref{GeneralLinearSystemODEs}, maps~\eqref{GeneralLinearSystemODEs} into the system
\[
\ddot{\boldsymbol{y}}=D\boldsymbol{y}, \qquad \mbox{where} \quad D=B-A^2.
\]
By~$J$ we denote the Jordan normal form of the matrix~$D$. Then there exists a nondegenerate matrix~$P$ such that $D=P^{-1}JP$,
and the point transformation $\boldsymbol{y}=P\boldsymbol{z}$
reduces the last system to the form $\ddot{\boldsymbol{z}}= J\boldsymbol{z}$.
As a result, for the study of symmetry properties of normal systems of linear second-order ODEs with constant coefficients it suffices
to consider only the systems of the form
\begin{equation}\label{basis_system}
\ddot{\boldsymbol{x}}=J\boldsymbol{x},
\end{equation}
where $J$ is a Jordan matrix,
\begin{gather}\label{Jordan_matrix}
J=\bigoplus_{l=1}^s J_{\lambda_l}^{k_l}, \qquad k_1+\cdots+k_s=n,
\end{gather}
$J_{\lambda_l}^{k_l}$ is the Jordan block of dimension~$k_l$ with the eigenvalue $\lambda_l$,
\[
\big[J_{\lambda_l}^{k_l}\big]_{ij}=\left\{
\begin{array}{ll}
\lambda_l, & \mbox{if}\ j=i, \\
1, & \mbox{if}\ j-i=1, \\
0, & \mbox{otherwise,}
\end{array}
\right. \quad i,j=1,\ldots,k_l,
\]
i.e., it is $k_l\times k_l$ matrix of the form
\[
J_{\lambda_l}^{k_l}=
\begin{pmatrix}
\lambda_l & 1 & 0 & 0 & \cdots & 0 \\
0 & \lambda_l & 1 & 0 & \cdots & 0 \\
0 & 0 & \lambda_l & 1 & \cdots & 0 \\
\cdots & \cdots & \cdots & \cdots & \cdots & \cdots\\
0 & 0 & 0 & 0 & \cdots & 1 \\
0 & 0 & 0 & 0 & \cdots & \lambda_l
\end{pmatrix},
\]
$(\lambda-\lambda_1)^{k_1}$, \dots, $(\lambda-\lambda_s)^{k_s}$ are elementary divisors of the matrix $J$.

In what follows, we also use the notation
$\diag(\gamma_1,\ldots,\gamma_l)$ for the $l\times l$ diagonal matrix
with the elements $\gamma_1$, \ldots, $\gamma_l$ on the diagonal,
and $E^l=\diag(1,\ldots,1)$ denotes the $l\times l$ unit matrix.
Subscripts of functions denote differentiation with respect to the corresponding variables.

\begin{remark}\label{note M=0}
If the matrix~$J$ is proportional to the unit matrix, $J\in \langle E^n\rangle$,
the associated system~\eqref{basis_system} is reduced to the form~\eqref{x''=0}
by a point transformation, which is called generalized Arnold's transformation~\cite{Gonzalez-Lopez1988}.
This is why we exclude such cases from the further consideration.
At the same time, the representation~\eqref{Jordan_matrix} of~$J$ as a direct sum of Jordan blocks
may contain partial direct sums of copies of the same $1\times1$ Jordan blocks,
and such partial direct sums are proportional to the unit matrices of the corresponding sizes.
\end{remark}

To compute Lie symmetries of system \eqref{basis_system}, we use the standard Lie approach \cite{Olver1993,Ovsiannikov1982}.
Acting by the second prolongation of the vector field
\[
Q=\xi(t,\boldsymbol{x})\p_t+\eta^a(t,\boldsymbol{x})\p_{x^a}
\]
on the system \eqref{basis_system} and substituting second-order derivatives using \eqref{basis_system}, we derive the
invariance condition
\begin{gather}
\eta^b_{tt}+2\eta^b_{x^at}x^a_t+\eta^b_{x^ax^c}x^a_tx^c_t+\eta^b_{x^a}(J\boldsymbol{x})^a-(\xi_{tt}+2\xi_{x^at}x^a_t+\xi_{x^ax^c}x^a_tx^c_t
+\xi_{x^a}(J\boldsymbol{x})^a)x^b_t\nonumber\\
\qquad{}-2(\xi_t+\xi_{x^a}x^a_t)(J\boldsymbol{x})^b=(J\boldsymbol{\eta})^b.\label{invariance_condition}
\end{gather}
Splitting \eqref{invariance_condition} with respect to the derivatives $x^a_t$ gives the system of determining equations
 for the coefficients $\xi(t,\boldsymbol{x})$ and $\eta^a(t,\boldsymbol{x})$,
\begin{gather}
 \xi_{x^ax^c}=0,
\label{3rd_degree}\\
 \eta^b_{x^ax^c}=0, \quad a\neq b\neq c, \qquad
 \eta^b_{x^ax^b}=\xi_{x^at}, \quad a\neq b, \qquad
 \eta^b_{x^bx^b}=2\xi_{x^b t},
\label{2nd_degree}\\
 \eta^b_{x^at}=\xi_{x^a} (J\boldsymbol{x})^b,\quad a\neq b, \qquad
 2\eta^b_{x^bt}=\xi_{tt}+2\xi_{x^b}(J\boldsymbol{x})^b+\xi_{x^a}(J\boldsymbol{x})^a,
\label{1st_degree}\\
 \eta^b_{tt}+\eta^b_{x^a}(J\boldsymbol{x})^a-2\xi_t(J\boldsymbol{x})^b=(J\boldsymbol{\eta})^b,
\label{0_degree}
\end{gather}
where there is no summation with respect to the repeated index~$b$.
Then equations~\eqref{3rd_degree} and~\eqref{2nd_degree} simultaneously imply
\begin{gather}\label{eq_xi+eta}
\xi=\xi^a(t)x^a+\xi^0(t), \qquad \eta^b=\xi^a_t(t)x^ax^b+\eta^{ba}(t)x^a+\eta^{b0}(t),
\end{gather}
where $\xi^a$, $\xi^0$, $\eta^{b a}$ and $\eta^{b0}$ are smooth functions of~$t$.
Substituting the expressions \eqref{eq_xi+eta} into the equations~\eqref{1st_degree} and~\eqref{0_degree}, and taking into the account the condition $J\notin \langle E^n\rangle$,
we obtain the more specific expressions
\[
\xi=c_1t+c_0, \qquad \eta^b=\eta^{ba}x^a+\eta^{b0}(t),
\]
where $c_1$, $c_0$, $\eta^{ba}$ are constants and $\eta^{b0}$ are smooth functions of~$t$.
The detail derivation of these expressions is the following.
The substitution of the expressions~\eqref{eq_xi+eta} into equations~\eqref{1st_degree}
and the subsequent partial splitting  with respect to~$\boldsymbol{x}$ result in the system
\begin{gather}\nonumber
\eta^{ba}_t=0, \quad
\xi^a_{tt}x^b=\xi^a(J\boldsymbol{x})^b, \quad a\ne b, \qquad
2\eta^{bb}_t=\xi^0_{tt}, \\ 
\label{2nd_degree_split4}
2\xi^b_{tt}x^b=-\xi^a_{tt}x^a+2\xi^b(J\boldsymbol{x})^b+\xi^a(J\boldsymbol{x})^a,
\end{gather}
where there is again no summation with respect to the repeated index~$b$.
Let $S=\{a_i,\, i=1,\dots,m\}$, $m\leqslant n$, be the subset of integers from $\{1,\dots, n\}$ for which $\xi^{a_i}\ne0$.
For any fixed $a_i\in S$ we have the condition
$\xi^{a_i}_{tt}/\xi^{a_i}=(J\boldsymbol{x})^b/x^b=:\mu_i=\mathrm{const}$ whenever $b\ne a_i$.
If $n>2$ and $m\geqslant2$, the constants $\mu_i$, $i=1,\dots,m$, coincide to each other.
If $n=m=2$, the above condition is equivalent to the set of equations
$(J\boldsymbol{x})^1=\mu_1x^1$, $(J\boldsymbol{x})^2=\mu_2x^2$, $\xi^1_{tt}=\mu_2\xi^1$ and $\xi^2_{tt}=\mu_1\xi^2$.
Then equations~\eqref{2nd_degree_split4} are reduced to the simple equations $(\mu_1-\mu_2)\xi^i=0$
and hence again $\mu_1=\mu_2$ as $\xi^1\xi^2\ne0$.
This is why both the cases considered contradict the condition $J\notin\langle E^n\rangle$.
If $n\geqslant2$ and $m=1$, only a single equation among~\eqref{2nd_degree_split4} is not satisfied identically,
which, after an additional arrangement, takes the form $(J\boldsymbol{x})^{a_1}=\mu_1x^{a_1}$ and implies,
jointly with the equations $(J\boldsymbol{x})^b=\mu_1x^b$ for $b\ne a_i$, the condition $J\in\langle E^n\rangle$
contradicting the supposition $J\notin\langle E^n\rangle$.
Therefore, the only possibility is $m=0$, i.e., $\xi^a=0$ for any~$a$.
Differentiating determining equations~\eqref{0_degree} with respect to~$t$ and then splitting  with respect to~$\boldsymbol{x}$,
we in particular obtain the matrix equation $\frac12\xi^0_{tttt}E^n-2\xi^0_{tt}J=0$
implying $\xi^0_{tt}=0$ in view of the linear independence of~$J$ and~$E^n$.
Then $\eta^{ba}_t=0$ for any pair of the indices $a$ and $b$.

Further, from the determining equations \eqref{0_degree} we also have
\begin{equation*}
\boldsymbol{\eta}^0_{tt}=J\boldsymbol{\eta}^0, \qquad \boldsymbol{\eta}^0=(\eta^{10},\dots,\eta^{n0})^{\rm T},
\end{equation*}
i.e.\ $\boldsymbol{\eta}^0$ is an arbitrary solution of the system \eqref{basis_system},
and additionally we have the following matrix equation for the $n\times n$ matrix $H=(\eta^{ba})$:
\begin{equation}\label{matrix equation}
HJ-2\xi_tJ=JH.
\end{equation}

If $\xi_t=0$ for any Lie symmetry operator of system~\eqref{basis_system},
equation~\eqref{matrix equation} is reduced to the condition of commutation of~$H$ with~$J$, $JH=HJ$.
Therefore, we arrive at the Frobenius problem: define all matrices $H$ that commute with a fixed matrix $J$.
This is a standard problem in matrix theory, see e.g.~\cite[Chapter~VIII]{Gantmaher}.

If there exists a Lie symmetry operator~$Q$ of system~\eqref{basis_system} with $\xi_t\not=0$,
the corresponding matrix~$H$ satisfies the inhomogeneous matrix equation~\eqref{matrix equation}
\begin{gather}\label{nilpotent-case1}
JH-HJ=\kappa J,
\end{gather}
where $\kappa=-2c_1\neq0$.
This equation is compatible if and only if the matrix~$J$ is nilpotent,
which directly follows from Lemma~4 of \cite[p.~44]{Jacobson} or from Theorem~II of~\cite{roth1952}.
The matrix equation~\eqref{nilpotent-case1} is an inhomogeneous linear system of algebraic equations
with respect to the coefficients of the matrix~$H$.
Hence the general solution of~\eqref{nilpotent-case1} is represented
as the sum of a particular solution of~\eqref{nilpotent-case1}
and the general solution of the associated homogeneous matrix equation $JH=HJ$
discussed in the previous paragraph.
A particular solution of the equation~\eqref{nilpotent-case1} is
$\kappa\diag(1,2,\dots,k_1,1,2,\dots,k_2,\,\dots,\,1,2,\dots,k_s)$,
where $k_1$, $k_2$, \dots, $k_s$ are the sizes of the Jordan blocks of the matrix~$J$, cf.~\eqref{Jordan_matrix}.

Summing up, we obtain the following theorem.

\begin{theorem}\label{main_theorem}
Suppose that $J$ is a matrix of the Jordan form~\eqref{Jordan_matrix} which is not proportional to the unit matrix.
The maximal Lie invariance algebra~$\mathfrak g^J$ of the system $\ddot{\boldsymbol{x}}=J\boldsymbol{x}$ is
\begin{gather*}
\langle \mathcal X^m, \, m=1,\dots,2n, \ \mathcal H^\ell, \, \ell=1,\dots, N, \ \mathcal T\rangle,
\end{gather*}
or
\begin{gather*}
\langle \mathcal X^m, \, m=1,\dots,2n, \ \mathcal H^\ell, \, \ell=1,\dots, N, \ \mathcal T, \ \mathcal D\rangle,
\end{gather*}
if $J$ is a non-nilpotent or nilpotent matrix, respectively.
Here
\begin{gather*}
\mathcal X^m= \varphi^{m a}(t)\partial_{x^a}, \qquad
\mathcal H^\ell= (H^\ell)^{ba}x^a\partial_{x^b},\qquad
\mathcal T=\partial_t, \qquad
\mathcal D=t\partial_t -2\gamma^{ab}x^b\partial_{x^a},
\end{gather*}
the vector-functions $\boldsymbol{\varphi}^m=(\varphi^{m 1}(t),\dots,\varphi^{m n}(t))^{\rm T}$, $m=1,\dots,2n$,
form a fundamental set of solutions for the system $\ddot{\boldsymbol{x}}=J\boldsymbol{x}$,
$H^\ell$, $\ell=1,\dots,N$, exhaust linearly independent matrices that commute with the matrix~$J$, and
$\gamma=(\gamma^{ab})=\diag (1,2,\dots,k_1,1,2,\dots,k_2,\,\dots,\,1,2,\dots,k_s)$.
\end{theorem}

By $N=N(D)$ we denote the number of linearly independent matrices that commute with an $n\times n$ matrix $D$.
It is obvious that $N(D)=N(\tilde D)$ if the matrices~$D$ and~$\tilde D$ are similar.

\begin{corollary}\label{dimension_corollary}
The dimension of the maximal Lie invariance algebra~$\mathfrak g^D$ of the system $\ddot{\boldsymbol{x}}=D\boldsymbol{x}$ with $D\notin \langle E^n\rangle$
is $2n+N+1$ or $2n+N+2$ if $D$ is a non-nilpotent or nilpotent matrix, respectively.
\end{corollary}

Let the matrix~$J$ of the form~\eqref{Jordan_matrix} be a Jordan form of the matrix~$D$ and
$\sigma_{ij}$ denote the degree of the greatest common divisor of the polynomials $(\lambda-\lambda_i)^{k_i}$ and $(\lambda-\lambda_j)^{k_j}$,
i.e., $\sigma_{ij}=0$ if $\lambda_i\ne\lambda_j$ and $\sigma_{ij}=\min(k_i,k_j)$ if $\lambda_i=\lambda_j$.
Then the number $N(D)$ is equal to \cite[p.~221]{Gantmaher}
\begin{gather}\label{dim_H}
N=\sum\limits_{i,j=1}^s \sigma_{ij}.
\end{gather}

Let $\mathcal I_1(\lambda)$, \dots, $\mathcal I_q(\lambda)$ form the complete set of nonconstant invariant polynomials of the matrix~$D$
with the degrees $n_1\geqslant \dots\geqslant n_q>0$.
Each invariant polynomial $\mathcal I_\alpha(\lambda)$ is a product of co-prime elementary divisors,
$\mathcal I_\alpha(\lambda)=(\lambda-\hat\lambda_1)^{d_{\alpha1}}\cdots(\lambda-\hat\lambda_p)^{d_{\alpha p}}$, $\alpha=1,\dots,q$.
Here
$\hat\lambda_1$, \dots, $\hat\lambda_p$ are all the distinct eigenvalues of the matrix $D$,
$d_{1j}\geqslant d_{2j}\geqslant\dots\geqslant d_{qj}\geqslant0$, $j=1,2,\dots,p$,
$n_\alpha=d_{\alpha1}+\dots+d_{\alpha p}$, $\alpha=1,\dots,q$,
$n_1+\dots+n_q=n$,
and, therefore,
$\mathcal I_1(\lambda)\cdots \mathcal I_q(\lambda)$ is the characteristic polynomial of the matrix $D$.
According to \cite[p.~222, Theorem~2]{Gantmaher}, we have one more representation for $N=N(D)$,
\begin{gather}\label{dim_H*}
N=n_1+3n_2 +\dots +(2q-1)n_q.
\end{gather}

We list other elementary properties of~$N(D)$, which are needed for the further consideration.
Thus, $N(D)=n\bmod2$, i.e.\ $N(D)$ takes only odd (resp.\ even) values for odd (resp.\ even)~$n$.
The value of~$N(D)$
is completely defined by the tuple $\bar n=(n_1,\dots,n_q)$ of the degrees of the nonconstant invariant polynomials of~$D$
or, equivalently, the partition of~$n$ into the integer summands $n_1\geqslant\dots\geqslant n_q>0$, $n=n_1+\dots +n_q$.
The representation~\eqref{dim_H*} implies that $n\leqslant N(D)\leqslant n^2$ for any $n\times n$ matrix~$D$.
The equality $N(D)=n$ holds if and only if $q=1$ and hence $\bar n=(n)$, i.e., all the elementary divisors of $D$ are co-prime in pairs
or, in other words, all the eigenvalues of~$D$ are distinct in pairs.
The maximal value $N=n^2$ is reached only if the matrix~$D$ is proportional to the unit matrix, $D\in \langle E^n\rangle$,
as then the number of (nonconstant) invariant polynomials associated with~$D$ is also maximal and equals~$n$, and $\bar n=(1,\dots,1)$.
The submaximal value of~$N(D)$ equals $N=n^2-2n+2$, and it is attained
only if $\bar n=(2,1,\dots,1)$; then $D$ is similar to either
$J_{\lambda_1}^2\oplus \big(\bigoplus_{i=1}^{n-2}J_{\lambda_1}^1\big)$ or
$\big(\bigoplus_{i=1}^{n-1}J_{\lambda_1}^1\big)\oplus J_{\lambda_2}^1$, where $\lambda_1\ne\lambda_2$.
The next less value $N=n^2-4n+8$ corresponds to the tuple $\bar n=(2,2,1,\dots,1)$, $n\geqslant 4$.

Theorem~\ref{main_theorem} and Corollary~\ref{dimension_corollary} jointly with~\eqref{dim_H} and~\eqref{dim_H*} give a much more effective algorithm
for computing the dimensions of the maximal Lie invariance algebras of systems from the class~\eqref{GeneralLinearSystemODEs}
than those existing in the literature, cf.~\cite[Proposition~4]{Campoamor2012}.
We also can explicitly construct a~basis of such algebra for a fixed Jordan matrix.
In particular, we directly derive series of simple estimates for dimensions of these algebras
(cf.\ also \cite{Campoamor2011,Campoamor2012}).

\begin{corollary}\label{CorollaryOnMinDimOfMIA}
The maximal Lie invariance algebra of the system $\ddot{\boldsymbol{x}}=D\boldsymbol{x}$ is of minimal dimension $3n + 1$ among systems of the form~\eqref{GeneralLinearSystemODEs}
if and only if the matrix~$D$ is not nilpotent and all the elementary divisors of~$D$ are co-prime in pairs.
\end{corollary}

In other words, the Jordan matrix associated with~$D$ consists of
either a single Jordan block with a nonzero eigenvalue
or a few Jordan blocks with eigenvalues different in pairs.
If $D$ is similar to the single Jordan block~$J_0^n$ with the zero eigenvalue,
the dimension of the maximal Lie invariance algebra of the system $\ddot{\boldsymbol{x}}=D\boldsymbol{x}$ equals $3n+2$.

\begin{corollary}
The dimensions of the maximal Lie invariance algebras of the systems of the form $\ddot{\boldsymbol{x}}=D\boldsymbol{x}$,
where $D\notin \langle E^n\rangle$, are not greater than $n^2+4$
and this upper bound is reached if and only if the matrix $D$ is similar to the nilpotent Jordan matrix
$J_0^2\oplus \big(\bigoplus_{i=1}^{n-2}J_0^1\big)$.
\end{corollary}

Therefore, for any system from the class~\eqref{GeneralLinearSystemODEs} which is not similar to the simplest system $\ddot{\boldsymbol{x}}=\boldsymbol{0}$
the dimension of its maximal Lie invariance algebra~$\mathfrak g$ satisfies the inequality
\[
3n+1\leqslant\dim\mathfrak g\leqslant n^2+4,
\]
and the bounds are exact.
The value $\dim\mathfrak g=n^2+3$ is attained for any system $\ddot{\boldsymbol{x}}=D\boldsymbol{x}$
with the matrix~$D$ similar to either
$J_{\lambda_1}^2\oplus \big(\bigoplus_{i=1}^{n-2}J_{\lambda_1}^1\big)$, where $\lambda_1\ne0$,
or
$\big(\bigoplus_{i=1}^{n-1}J_{\lambda_1}^1\big)\oplus J_{\lambda_2}^1$, where $\lambda_1\ne\lambda_2$,
and only for elements from the equivalence subclasses of such systems with respect to point transformations.
As above mentioned, we have $\dim\mathfrak g=3n+2$ if $D$ is similar to $J_0^n$.
At the same time, only in the case $2\leqslant n\leqslant4$ for each integer value~$\rho$ from the interval $[3n+1,n^2+4]$
there exists a system from the class~\eqref{GeneralLinearSystemODEs} whose maximal Lie invariance algebra is of dimension~$\rho$.
For $n\geqslant 5$, there exists, in particular, no $n\times n$ matrix~$D$ with $\dim \mathfrak g^D\in [n^2-2n+11,n^2+2]$.
The number of such exceptional intervals grows with increasing~$n$.

If the matrix $J$ is diagonal, i.e.\ all its elementary divisors are of degree~$1$, and additionally $J\notin \langle E^n\rangle$ then
we have $N=N(J)=\sum_{i=1}^p r_i^2$,
where $\hat\lambda_1$, \dots, $\hat\lambda_p$ are all the distinct eigenvalues of the matrix $J$ and
$r_i$ is the multiplicity of~$\smash{\hat\lambda_i}$, $i=1,\dots,p$.
This is why Corollary~\ref{dimension_corollary} directly implies Proposition~3 of~\cite{Campoamor2011}.

{\samepage Given a Jordan matrix $J$ with a single eigenvalue $\lambda$
(i.e. $J=J_\lambda^{k_1}\oplus J_\lambda^{k_2}\oplus\cdots \oplus J_\lambda^{k_s},$
where $k_1\geqslant k_2\geqslant\cdots \geqslant k_s$ and $k_1+\dots +k_s=n$),
it follows from~\eqref{dim_H} and~\eqref{dim_H*} that
\[
N=\sum_{i=1}^{s}(2i-1)k_i=ns-\sum_{i=1}^{s-1} \sum_{j=i+1}^s (k_i-k_j).
\]
In view of Corollary~\ref{dimension_corollary}, this essentially simplifies Theorem~2 of~\cite{Campoamor2012}.

}

\begin{remark}\label{RemarkOnRealCase}
Given a matrix~$D$ with real entries,
the number of linearly independent solutions of the systems $\ddot{\boldsymbol{x}}=D\boldsymbol{x}$
(resp.\ the number $N=N(D)$ of linearly independent matrices commuting with~$D$) over the real field is the same as over the complex field.
Therefore, all results of this section are directly extended to the real case.
Thus, in Theorem~\ref{main_theorem} it suffices to consider real Jordan matrices and to take real counterparts for
the vector-functions $\boldsymbol{\varphi}^m=(\varphi^{m 1}(t),\dots,\varphi^{m n}(t))^{\rm T}$, $m=1,\dots,2n$,
and the matrices $H^\ell$, $\ell=1,\dots,N$.
See the second example in the next section.
\end{remark}

\section{Illustrative examples}\label{examples}

We present two simple examples that illustrate Theorem~\ref{main_theorem}.

\begin{example}
For system~\eqref{basis_system} with the Jordan matrix $J=J_{\lambda_1}^2\oplus J_{\lambda_2}^2$, i.e., the system
\begin{gather}\label{example1}
\ddot x^1=\lambda_1x^1+x^2,\quad
\ddot x^2=\lambda_1x^2,\quad
\ddot x^3=\lambda_2x^3+x^4,\quad
\ddot x^4=\lambda_2x^4,
\end{gather}
there are two different cases subject to the eigenvalues $\lambda_1$ and~$\lambda_2$, namely $\lambda_1\not=\lambda_2$ and $\lambda_1=\lambda_2$.

For $\lambda_1\not=\lambda_2$ the single nonconstant invariant polynomial of the matrix $J$ is $(\lambda-\lambda_1)^2 (\lambda-\lambda_2)^2$.
Hence $N=4$, and any matrix commuting with~$J$ has the form \cite[Chapter~VIII]{Gantmaher}
\[
  H=
\begin{pmatrix}
\eta^{11}&\eta^{12}& \!\!\!\vdots\!\!\! & 0&0 \\[-1ex]
0& \eta^{11} &\!\!\!\vdots\!\!\! & 0 & 0 \\[-1ex]
\hdotsfor{5}\\[-1ex]
0 & 0 & \!\!\!\vdots\!\!\! &\eta^{33} & \eta^{34} \\[-1ex]
0 & 0 & \!\!\!\vdots\!\!\! & 0 & \eta^{33}
\end{pmatrix}.
\]
Finding a fundamental set of solutions of the system \eqref{example1} and using Theorem~\ref{main_theorem},
we obtain the 13-dimensional maximal Lie invariance algebras of the system \eqref{example1}
\begin{gather*}
 \langle e^{\sqrt{\lambda_1}t}(t\p_{x^1}+2\sqrt{\lambda_1}\p_{x^2}), \ e^{-\sqrt{\lambda_1}t}(t\p_{x^1}-2\sqrt{\lambda_1}\p_{x^2}), \ e^{\sqrt{\lambda_1}t}\p_{x^1}, \ e^{-\sqrt{\lambda_1}t}\p_{x^1}, \\
 e^{\sqrt{\lambda_2}t}(t\p_{x^3}+2\sqrt{\lambda_2}\p_{x^4}), \ e^{-\sqrt{\lambda_2}t}(t\p_{x^3}-2\sqrt{\lambda_2}\p_{x^4}), \ e^{\sqrt{\lambda_2}t}\p_{x^3}, \ e^{-\sqrt{\lambda_2}t}\p_{x^3}, \\
  x^1\p_{x^1}+x^2\p_{x^2}, \ x^2\p_{x^1}, \ x^3\p_{x^3}+x^4\p_{x^4}, \ x^4\p_{x^3}, \ \p_t
 \rangle
 \end{gather*}
for $\lambda_1\neq\lambda_2\neq0$ or {\samepage
\begin{gather*}
 \langle e^{\sqrt{\lambda_1}t}(t\p_{x^1}+2\sqrt{\lambda_1}\p_{x^2}), \ e^{-\sqrt{\lambda_1}t}(t\p_{x^1}-2\sqrt{\lambda_1}\p_{x^2}), \ e^{\sqrt{\lambda_1}t}\p_{x^1}, \ e^{-\sqrt{\lambda_1}t}\p_{x^1}, \ t^3\p_{x^3}+6t\p_{x^4}, \\
 t^2\p_{x^3}+2\p_{x^4}, \ t\p_{x^3}, \ \p_{x^3}, \ x^1\p_{x^1}+x^2\p_{x^2}, \ x^2\p_{x^1}, \ x^3\p_{x^3}+x^4\p_{x^4}, \ x^4\p_{x^3}, \ \p_t
\rangle
\end{gather*}
for $\lambda_1\neq\lambda_2= 0$.}

If $\lambda_1=\lambda_2$, the nonconstant invariant polynomials of the matrix $J$ are $(\lambda-\lambda_1)^2$ and $(\lambda-\lambda_1)^2$,
and hence $N=8$. A matrix~$H$ commutes with~$J$ if and only if it has the~form
\[
H=
\begin{pmatrix}
\eta^{11}&\eta^{12}& \!\!\!\vdots\!\!\!& \eta^{13}&\eta^{14} \\[-1ex]
0& \eta^{11}& \!\!\!\vdots\!\!\! & 0 &\eta^{13} \\[-1ex]
\hdotsfor{5}\\[-1ex]
\eta^{31} & \eta^{32}& \!\!\!\vdots\!\!\! &\eta^{33} & \eta^{34} \\[-1ex]
0 & \eta^{31}& \!\!\!\vdots\!\!\! & 0 & \eta^{33}
\end{pmatrix}.
\]
As a result, we construct the 17-dimensional maximal Lie invariance algebra
\begin{gather*}
\langle
e^{\sqrt{\lambda_1}t}(t\p_{x^1}+2\sqrt{\lambda_1}\p_{x^2}), \ e^{-\sqrt{\lambda_1}t}(t\p_{x^1}-2\sqrt{\lambda_1}\p_{x^2}), \ e^{\sqrt{\lambda_1}t}\p_{x^1}, \ e^{-\sqrt{\lambda_1}t}\p_{x^1}, \\
e^{\sqrt{\lambda_1}t}(t\p_{x^3}+2\sqrt{\lambda_1}\p_{x^4}), \ e^{-\sqrt{\lambda_1}t}(t\p_{x^3}-2\sqrt{\lambda_1}\p_{x^4}), \ e^{\sqrt{\lambda_1}t}\p_{x^3}, \ e^{-\sqrt{\lambda_1}t}\p_{x^3}, \ x^1\p_{x^1}+x^2\p_{x^2}, \\
x^2\p_{x^1}, \ x^3\p_{x^1}+x^4\p_{x^2}, \ x^4\p_{x^1}, \ x^1\p_{x^3}+x^2\p_{x^4}, \ x^2\p_{x^3}, \ x^3\p_{x^3}+x^4\p_{x^4}, \ x^4\p_{x^3}, \ \p_t
\rangle
\end{gather*}
if $\lambda_1\neq 0$ or the 18-dimensional maximal Lie invariance algebra
\begin{gather*}
\langle
t^3\p_{x^1}+6t\p_{x^2}, \ t^2\p_{x^1}+2\p_{x^2}, \ t\p_{x^1}, \ \p_{x^1}, \ t^3\p_{x^3}+6t\p_{x^4}, \ t^2\p_{x^3}+2\p_{x^4}, \ t\p_{x^3}, \ \p_{x^3}, \\
x^1\p_{x^1}+x^2\p_{x^2}, \ x^2\p_{x^1}, \ x^3\p_{x^1}+x^4\p_{x^2}, \ x^4\p_{x^1}, \ x^1\p_{x^3}+x^2\p_{x^4}, \ x^2\p_{x^3}, \ x^3\p_{x^3}+x^4\p_{x^4},  \\
x^4\p_{x^3}, \ \p_t, \ t\p_t-2x^1\p_{x^1}-4x^2\p_{x^2}-2x^3\p_{x3}-4x^4\p_{x^4}
\rangle
\end{gather*}
if $\lambda_1 = 0$.
\end{example}

\begin{example}
Over the real field, consider the system
\begin{gather}\label{example2}
\ddot x_1=\mu x_1+\nu x_2,\quad
\ddot x_2=-\nu x_1+\mu x_2,\quad
\ddot x_3=x_4,\quad
\ddot x_4=0,\quad
\ddot x_5=0,
\end{gather}
which is of the form~\eqref{basis_system} with the real-valued Jordan matrix
$J=R_{\mu\nu}^2 \oplus J_{0}^2 \oplus J_{0}^1$, where
$R_{\mu\nu}^2=\begin{pmatrix}
\mu & \nu \\
-\nu & \mu
\end{pmatrix}$.
A matrix~$H$ commutes with~$J$ if and only if it has the~form
\[
H=
\begin{pmatrix}
\eta^{11} & \eta^{12} & \!\!\!\vdots\!\!\! & 0 & 0 & \!\!\!\vdots\!\!\! & 0\\
-\eta^{12} & \eta^{11} & \!\!\!\vdots\!\!\! & 0 &0 & \!\!\!\vdots\!\!\! &0 \\[-1ex]
\hdotsfor{7}\\[-1ex]
0  & 0 & \!\!\!\vdots\!\!\! &\eta^{33} & \eta^{34} & \!\!\!\vdots\!\!\! & \eta^{35}\\
0  & 0 & \!\!\!\vdots\!\!\! & 0 & \eta^{33}& \!\!\!\vdots\!\!\! &0 \\[-1ex]
\hdotsfor{7}\\[-1ex]
0 & 0& \!\!\!\vdots\!\!\! & 0 & \eta^{54} & \!\!\!\vdots\!\!\! & \eta^{55}
\end{pmatrix}.
\]
Then Theorem~\ref{main_theorem} implies
that the maximal Lie invariance algebra of the system \eqref{example2} is generated by the 18 vector fields
\begin{gather*}
\nu e^{\alpha t} \cos \beta t \p_{x^1}+e^{\alpha t} \big((\alpha^2-\beta^2-\mu)\cos\beta t -2\alpha\beta\sin\beta t\big) \p_{x^2} , \\
\nu e^{\alpha t} \sin \beta t \p_{x^1}+e^{\alpha t} \big((\alpha^2-\beta^2-\mu)\sin\beta t + 2\alpha\beta\cos\beta t\big) \p_{x^2} ,\\
\nu e^{-\alpha t} \cos \beta t \p_{x^1}+e^{-\alpha t} \big((\alpha^2-\beta^2-\mu)\cos\beta t + 2\alpha\beta\sin\beta t\big) \p_{x^2} , \\
\nu e^{-\alpha t} \sin \beta t \p_{x^1}+e^{-\alpha t} \big((\alpha^2-\beta^2-\mu)\sin\beta t - 2\alpha\beta\cos\beta t\big) \p_{x^2} ,\\
 \p_{x^3}, \ t\p_{x^3}, \ t^3\p_{x^3}+6t\p_{x^4}, \ t^2\p_{x^3}+2\p_{x^4}, \ \p_{x^5}, \ t\p_{x^5}, \\
 x^1\p_{x^1}+x^2\p_{x^2}, \ x^2\p_{x^1} - x^1\p_{x^2}, \ x^3\p_{x^3}+x^4\p_{x^4}, \ x^4\p_{x^3}, \ x^5\p_{x^3}, \ x^4\p_{x^5}, \ x^5\p_{x^5}, \ \p_t,
\end{gather*}
where $\alpha=(\mu^2+\nu^2)^{1/4}\cos \big(\frac 12 \arctan\frac{\nu}{\mu}\big)$, $\beta=(\mu^2+\nu^2)^{1/4}\sin \big(\frac 12 \arctan\frac{\nu}{\mu}\big)$.
\end{example}

\section{Conclusion}

In this paper we exhaustively study Lie symmetries of systems of second-order linear ordinary differential equations with commuting constant matrix coefficients
over both the complex and real fields, which reduce to the form~\eqref{basis_system}.
The explicit description of the maximal Lie invariance algebra
of any system from the class~\eqref{GeneralLinearSystemODEs} is presented by Theorem~\ref{main_theorem}.
Corollary~\ref{dimension_corollary} jointly with formulas~\eqref{dim_H} and~\eqref{dim_H*}
give a simple and algorithmic tool in order for computing the dimensions of such algebras.
In particular, we show that these dimensions are completely defined by the degrees of the nonconstant invariant polynomials
of the corresponding Jordan matrices in the reduced form~\eqref{basis_system}.
We also give estimates for possible values of these dimensions.
In order to fix ideas, the consideration is carried out over the complex field
but it can be directly extended to the real case with preserving all assertions, cf.\ Remark~\ref{RemarkOnRealCase}.
The results obtained in the present paper amend and generalize those of \cite{Campoamor2011,Campoamor2012,Meleshko2011,WafoSoh2010}.
The advantages and simplicity of the approach proposed are illustrated by examples.

As a next step, we plan to study Lie symmetries of more general systems than \eqref{GeneralLinearSystemODEs},
in particular, classes of linear systems with variable coefficients.
This will need the investigation of the equivalence groupoids of the above classes, the description of their maximal normalized subclasses and
combining the algebraic and compatibility methods of group classification~\cite{Popovych2010}.
A~still open problem is also the linearization of systems \eqref{second_order_system} by point transformations.

\subsection*{Acknowledgments}

We thank Professors Yurii Samoilenko and Vasyl Ostrovskyi for helpful discussions on matrix equations.
The research of ROP was supported by the Austrian Science Fund (FWF), project P20632.
It is our great pleasure to thank the referee for useful suggestions that have considerably improved the paper.

\end{document}